\documentclass[letterpaper, 10 pt, conference]{ieeeconf}  % Comment this line out if you need a4paper

\IEEEoverridecommandlockouts                              % This command is only needed if
\usepackage{graphics} % for pdf, bitmapped graphics files
\usepackage{graphicx}
\usepackage{epsfig} % for postscript graphics files
\usepackage{mathptmx} % assumes new font selection scheme installed
\usepackage{times} % assumes new font selection scheme installed
\usepackage{amsmath,amssymb,amsfonts}
\usepackage{url}
\usepackage{algorithm,algorithmicx}
\usepackage{algpseudocode}
\usepackage{lineno}
\usepackage{hyperref}
\usepackage{bbm}
\usepackage{cite}
\title{{\LARGE{ \bf
Solving Markov decision processes for network-level post-hazard recovery via simulation optimization and rollout}}\\ \normalsize\emph{(Invited Paper)}
}

\author{Yugandhar Sarkale$^{1}$, Saeed Nozhati$^{2}$, Edwin K. P. Chong$^{1}$, Bruce R. Ellingwood$^{2}$, Hussam Mahmoud$^{2}$% <-this % stops a space
\thanks{$^{1}$Yugandhar Sarkale and Edwin K. P. Chong are with Department of Electrical and Computer Engineering, Colorado State University, Fort Collins, CO 80523-1373, USA {\tt\small Yugandhar.Sarkale, Edwin.Chong@colostate.edu}}%
\thanks{$^{2}$Saeed Nozhati, Bruce R. Ellingwood and Hussam Mahmoud are with the Department of Civil and Environmental Engineering, Colorado State University, Fort Collins, CO 80523-1372, USA {\tt\small Saeed.Nozhati, Bruce.Ellingwood, Hussam.Mahmoud@colostate.edu}}%
\thanks{This work was supported by the National Science Foundation under Grant CMMI-1638284. This support is gratefully acknowledged. Any opinions, findings, conclusions, or recommendations presented in this material are solely those of the authors and do not necessarily reflect the views of the National Science Foundation.}%
}
\begin{document}

\maketitle
\thispagestyle{empty}
\pagestyle{empty}

%%%%%%%%%%%%%%%%%%%%%%%%%%%%%%%%%%%%%%%%%%%%%%%%%%%%%%%%%%%%%%%%%%%%%%%%%%%%%%%%
\begin{abstract}
 Computation of optimal recovery decisions for community resilience assurance post-hazard is a combinatorial decision-making problem under uncertainty. It involves solving a large-scale optimization problem, which is significantly aggravated by the introduction of uncertainty. In this paper, we draw upon established tools from multiple research communities to provide an effective solution to this challenging problem. We provide a stochastic model of damage to the water network (WN) within a testbed community following a severe earthquake and compute near-optimal recovery actions for restoration of the water network.  We formulate this stochastic decision-making problem as a Markov Decision Process (MDP), and solve it using a popular class of heuristic algorithms known as \emph{rollout}. A simulation-based representation of MDPs is utilized in conjunction with rollout and the Optimal Computing Budget Allocation (OCBA) algorithm to address the resulting stochastic simulation optimization problem. Our method employs non-myopic planning with efficient use of simulation budget. We show, through simulation results, that rollout fused with OCBA performs competitively with respect to rollout with total equal allocation (TEA) at a meagre simulation budget of 5-10\% of rollout with TEA, which is a crucial step towards addressing large-scale community recovery problems following natural disasters.

\end{abstract}

%%%%%%%%%%%%%%%%%%%%%%%%%%%%%%%%%%%%%%%%%%%%%%%%%%%%%%%%%%%%%%%%%%%%%%%%%%%%%%%%
\section{INTRODUCTION}
Natural disasters have a significant impact on the economic, social, and cultural fabric of affected communities. Moreover, because of the interconnected nature of communities in the modern world, the adverse impact is no longer restricted to the locally affected region, but it has ramifications on national or international scale. Among other factors, the occurrence of such natural disasters is on the rise owing to population growth and economic development in hazard-prone areas \cite{saeed}. Keeping in view the increased frequency of natural disasters, there is an urgent need to address the problem of community recovery post-hazard. Typically, the resources available to post-disaster planners are limited and relatively small compared to the impact of the damage. Under these scenarios, it becomes imperative to assign limited resources to various damaged components in the network optimally to support community recovery. Such an assignment must also consider multiple objectives and cascading effects due to the interconnectedness of various networks within the community and must also successfully adopt previous proven methods and practices developed by expert disaster-management planners. Holistic approaches addressing various uncertainties for network-level management of limited resources must be developed for maximum effect. Civil infrastructure systems, including power, transportation, and water networks, play a critical part in post-disaster recovery management. In this study, we focus on one such critical infrastructure system, namely the water networks (WN), and compute near-optimal recovery actions, in the aftermath of an earthquake, for the WN of a test-bed community.

Markov decision processes (MDPs) offer a convenient framework for representation and solution of stochastic decision-making problems. Exact solutions are intractable for problems of even modest size; therefore, approximate solution methods have to be employed. We can leverage the rich theory of MDPs to model recovery action optimization for large state-space decision-making problems such as our. In this study, we employ a simulation-based representation and solution of MDP. The near-optimal solutions are computed using an approximate solution technique known as \emph{rollout}. Even though state-of-the-art hardware and software practices are used to implement the solution to our problem, we are faced with the additional dilemma of computing recovery actions on a fixed simulation budget without affecting the solution performance. Therefore, any prospective methodology must incorporate such a limitation in its solution process. We incorporate the Optimal Computing Budget Allocation (OCBA) algorithm into our MDP solution process \cite{ocbar1,ocbarr} to address the limited simulation budget problem.

\section{TESTBED CASE STUDY}

\subsection{Network Characterization}\label{test1}
\begin{figure}[t]
	
	\centering
	\includegraphics[width=\columnwidth]{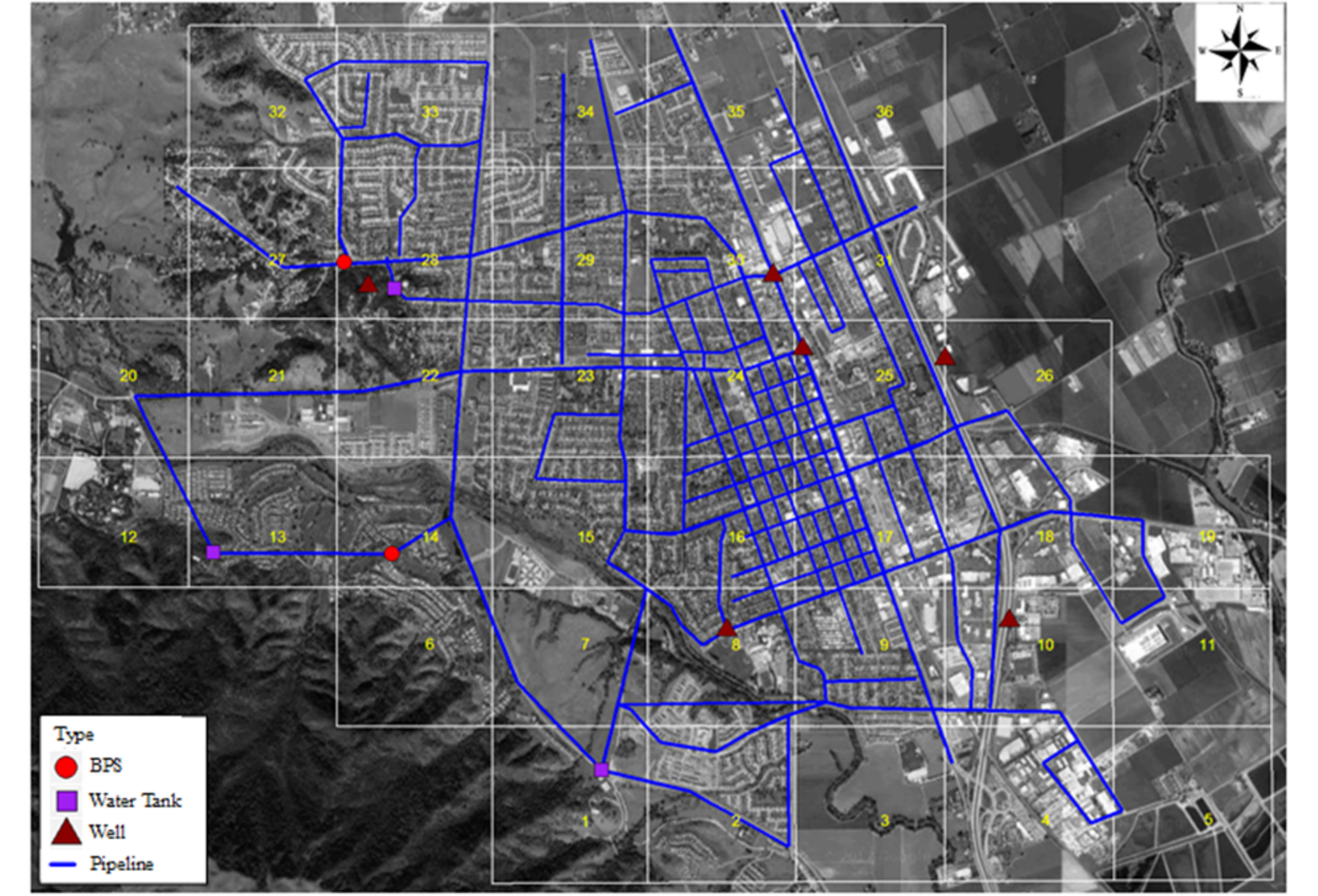}
	\caption{The modeled Water Network of Gilroy}
	\label{fig1}
\end{figure}

This study considers the potable water network (WN) of Gilroy, CA, USA as an example to illustrate the proposed methodology. Gilroy, located 50 kilometers (km) south of the city of San Jose, CA is approximately 41.91 km\textsuperscript{2} in area, with a population of 48,821\cite{Gilroy1}. We divide our study area into 36 grid regions to define the properties of infrastructure systems, household units, and the population. Our model of the community maintains adequate detail to study the performance of the WN at a community level under severe earthquakes. The potable water of Gilroy is provided only by the Llgas sub-basin \cite{Gilroy2}. The potable water wells, located in wood-frame buildings, pump water into the distribution system. The Gilroy municipal water pipelines range from 102~mm to 610~mm in diameter \cite{Gilroy2}. In this study, a simplified aggregated model of WN of Gilroy adopted from \cite{Gilroy2} is modeled. This model shown in Fig.~\ref{fig1}, includes six water wells, two booster pump stations (BPS), three water tanks (WT), and the main pipelines.
%%%%%%%%%%%%%%%%%%%%%%%%%%%   Fig 1  %%%%%%%%%%%%%%%%%%%%%%
%%%%%%%%%%%%%%%%%%%%%%%%%%%%%%%%%%%%%%%%%%%%%%%%%%%%%%%%%%
%%%%%%%%%%%%%%%%%%%%%%%%%%%   Fig 2  %%%%%%%%%%%%%%%%%%%%%%
%\begin{figure}[t]
	
%	\centering
%	\caption{The simulation of median of peak ground acceleration field}
%	\label{fig2}
%\end{figure}
%%%%%%%%%%%%%%%%%%%%%%%%%%%%%%%%%%%%%%%%%%%%%%%%%%%%%%%%%%
%%%%%%%%%%%%%%%%%%%%%%%%%%%   Fig 3  %%%%%%%%%%%%%%%%%%%%%%
%\begin{figure}[t]
	
%	\centering
%	\includegraphics[width=\linewidth]{PGV}
%	\caption{The simulation of median of peak ground velocity field}
%	\label{fig3}
%\end{figure}
%%%%%%%%%%%%%%%%%%%%%%%%%%%%%%%%%%%%%%%%%%%%%%%%%%%%%%%%%%

\subsection{Seismic Hazard Simulation}\label{test2}
The San Andreas Fault (SAF), which is near Gilroy, is a source of severe earthquakes. In this study, we assume that a seismic event of moment magnitude $Mw=6.9$ occurs at one of the closest points on the SAF projection to downtown Gilroy with an epicentral distance of approximately 12 km. Ground motion prediction equations (GMPE) determine the conditional probability of exceeding ground motion intensity at specific geographic locations within Gilroy for this earthquake.

The Abrahamson et al. \cite{abrahamson} GMPE is used to estimate the Intensity Measures (IM) and associated uncertainties. Peak Ground Acceleration (PGA) is considered for the above-ground WN facilities and wells, whereas Peak Ground Velocity (PGV) is considered as IM of pipelines.

\subsection{Fragility and Restoration Assessment of Water Network}\label{test3}
The physical damage to WN components can be assessed by seismic fragility curves. We use the fragility curves presented in HAZUS-MH \cite{hazus} for wells, water tanks, and pump stations based on the IM of PGA. This study adopts the assumptions in \cite{adachi} for water pipelines. The failure probability of a pipe is bounded as follows:

\begin{equation}\label{eq1}
1-G_{\varepsilon PGV}(-CL\mu_{PGV})\leq E[P_{f}] \leq 1-E[exp(-CL\mu_{PGV})]
\end{equation}
where $P_{f}$ is the failure probability of a pipe, $L$ is the length of pipe, $\mu_{PGV}$ is the average PGV for the entire length of the water main, and $G(\cdot)$ is the moment-generating function of $\varepsilon PGV$ (the residual of the PGV). The term $C$ for water pipe segment \textit{i} is $C=K\times0.00187\times PGV_{\textit{i}}$, where $K$ is a coefficient determined by the pipe material, diameter, joint type, and soil condition based on the guidelines prepared by the American Lifeline Alliance \cite{ALA}. Adachi and Ellingwood \cite{adachi} demonstrated that the Upper Bound (UB) and exact solutions \eqref{eq1} are close enough so that in practical applications the UB assessment (conservative evaluation) can be used.

Repair crews, replacement components, and tools are considered as available units of resources to restore the damaged components of WN following the hazard. One unit of resources is required to repair each damaged component \cite{ouyang,masoomi}. However, the available units of resources are limited and depend on the capacities and policy of the entities within the community. To restore the WN, the restoration times based on exponential distributions synthesized from HAZUS-MH \cite{hazus} are used, as summarized in Table~\ref{T1}. The pipe-restoration time in the WN is based on repair rate or number of repairs per kilometer.

%%%%%%%%%%%%%%%%%   Table 1 %%%%%%%%%%
\begin{table}[h]
	\caption{The expected repair times (Unit: days)}\label{T1}
	\resizebox{\linewidth}{!}{
		\begin{tabular}{lllll}
			\hline
			&Damage States\\
			\hline
			Component & Minor & Moderate & Extensive & Complete \\
			\hline
			Water tanks & 1.2 & 3.1 & 93 & 155 \\
			Wells & 0.8 & 1.5 & 10.5 & 26 \\
			Pumping plants & 0.9 & 3.1 & 13.5 & 35 \\
			\hline
		\end{tabular}
	}
\end{table}
%%%%%%%%%%%%%%%%%%%%%%%%%%%%%%%

\section{PROBLEM DESCRIPTION and SOLUTION}
\subsection{MDP Framework}
We provide a brief description of MDP \cite{Puterman} for the sake of completeness. An MDP is a controlled dynamical process useful in modelling of wide range of decision-making problems. It can be represented by the 4-tuple $\langle S,A,T,R \rangle$. Here, $S$ represents the set of states, and $A$ represents the set of actions. Let $s,s' \in S$ and $a \in A$; then $T$ is the state transition function, where $T(s,a,s')=P(s'\mid s,a)$ is the probability of going into state $s'$ after taking action $a$ in state $s$. $R$ is the reward function, where $R(s,a,s')$ is the reward received after transitioning from $s$ to $s'$ as a result of action $a$. In this study, we assume that $|S|$ and $|A|$ are finite; $R$ is bounded and real-valued and a deterministic function of $s$, $a$ and $s'$. Implicit in our presentation are also the following assumptions: First order Markovian dynamics (history independence), stationary dynamics (reward function is not a function of absolute time), and full observability of the state space (outcome of an action in a state might be random, but we know the state reached after action is completed). In our study, we assume that we are allowed to take recovery actions (decisions) indefinitely until all the damaged components of our modeled problem are repaired (infinite-horizon planning). In this setting, we have a stationary policy $\pi$, which is defined as $\pi: S \rightarrow A$. Suppose that decisions are made at discrete-time $t$; then $\pi(s)$ is the action to be taken in state $s$ (regardless of time $t$). Our objective is to find an optimal policy $\pi^*$. For the infinite-horizon case, $\pi^*$ is defined as
\begin{equation}\label{opt}
  \pi^*=\arg\max_{\pi} V^{\pi}(s_0),
\end{equation}
where
\begin{equation}\label{val}
 V^\pi(s_0)=E\left\lbrack\sum_{t=0}^{\infty}\gamma^{\,t}R(s_t,\pi(s_t),s_{t+1})\right\rbrack
\end{equation}
is called the value function for a fixed policy $\pi$, and $0<\gamma<1$ is the discount factor. Note that the optimal policy is independent of the initial state $s_0$. Also, note that we maximize over policies $\pi$, where at each time $t$ the action taken is $a_t=\pi(s_t)$. Stationary optimal policies are guaranteed to exist for discounted infinite-horizon optimization criteria \cite{howard}. To summarize, our presentation is for infinite-horizon discrete-time MDPs with the discounted value as our optimization criterion.

\subsection{MDP Solution}\label{sol}
A solution to an MDP is the optimal policy $\pi^*$. We can obtain $\pi^*$ with linear programming or dynamic programming. In the dynamic programming regime, there are several solution strategies, namely value iteration, policy iteration, modified policy iteration, etc. Unfortunately, such exact solution algorithms are intractable for large state and actions spaces. We briefly mention here the method of value iteration because it illustrates the Bellman's equation \cite{Bellman}. Studying Bellman's equation is useful for defining $Q$ value function. $Q$ value function will play a critical role in describing the rollout algorithm. Let $V^{\pi^{*}}$ denote the optimal value function for some $\pi^*$; Bellman showed that $V^{\pi^{*}}$ satisfies:
\begin{align}\label{bell}
  V^{\pi^{*}}(s)= \max_{a \in A(s)}\left \{\gamma \cdot \sum_{s'}P(s'\mid s,a)\cdot \left \lbrack V^{\pi^{*}}(s') +R(s,a,s')\right \rbrack\right \}.
\end{align}
Equation \eqref{bell} is known as the Bellman's optimality equation, where $A(s)$ is the set of possible actions in any state $s$. The value iteration algorithm solves \eqref{bell} by using Bellman backup repeatedly, where Bellman backup is given by:
\begin{equation}\label{bellu}
  V_{i+1}(s)=\max_{a \in A(s)}\left\{\gamma \sum_{s'} P(s'\mid s,a)\cdot \left \lbrack V_{i}(s')+R(s,a,s')\right \rbrack \right\}.
\end{equation}
Bellman showed that $\lim_{i\to\infty} V_i = V^{\pi^*}$, where $V_0$ is initialised arbitrarily.\footnote{On a historical note, Lloyd Shapely's paper \cite{Shapley} included the value iteration algorithm for MDPs as a special case, but this was recognised only later on \cite{kallenberg2003finite}.} Next, we define the $Q$ value function of policy $\pi$:
\begin{equation}\label{Qval}
  Q_{\pi}(s,a)=\gamma \cdot \sum_{s'}P(s'\mid s,a)\cdot \left \lbrack V^{\pi}(s') + R(s,a,s')\right \rbrack.
\end{equation}
 The $Q$ value function of any policy $\pi$ gives the expected discount reward in the future after starting in some state $s$, taking action $a$ and following policy $\pi$ thereafter. Note that this is the inner term in \eqref{bell}.

\subsection{Simulation-Based Representation of MDP}
We now briefly explain the simulation-based representation of an MDP \cite{Fern}. Such a representation serves well for large state and action spaces, which is a characteristic feature of many real-world problems. When $|S|$ or $|A|$ is large, it is not feasible to represent $T$ and $R$ in a matrix form. A simulation-based representation of an MDP is a 5-tuple $\langle S,A,R,T,I \rangle$, where $S$ and $A$ are as before, except $|S|$ and $|A|$ are large. Here, $R$ is a stochastic real-valued bounded function that stochastically returns a reward $r$ when input $s$ and $a$ are provided, where $a$ is the action applied in state $s$. $T$ is a simulator that stochastically returns a state $s'$ when state $s$ and action $a$ are provided as inputs. $I$ is the stochastic initial state function that stochastically returns a state according to some initial state distribution. $R$, $T$, and $I$ can be thought of as any callable library functions that can be implemented in any programming language.

\subsection{Problem Formulation}\label{probform}
After an earthquake event occurs, the components of the water network remain undamaged or exhibit one of the damage states as shown in Table~\ref{T1}. Let $L'$ be the total number of damaged component at $t$. Let $t_{c}$ represent the decision time when all components are repaired. There is a fixed number of resource units ($M$) available to the decision maker. At each discrete-time $t$, the decision maker has to decide the assignment of unit of resource to the damaged locations; each component cannot be assigned more than one resource unit. When the number of damaged locations is less than the number of units of resources (because of sequential application of repair actions, or otherwise), we retire the extra unit of resources so that $M$ is equal to the number of damaged locations.
\begin{itemize}
  \item \emph{States S:} Let $s_t$ be the state of the damaged components of the system at time $t$; then $s_t$ is a vector of length $L'$, $s_t=(s_t^1,\ldots,s_t^{L'})$, and $s_t^l$ is one of the damaged state in Table~\ref{T1} where  $l\in \{1,\ldots,L' \}$.
  \item \emph{Actions A:} Let $a_t$ denote the repair action to be carried out at time $t$. Then, $a_t$ is a vector of length $L'$, $a_t=(a_t^1, \ldots, a_t^{L'})$, and $a_t^l \in \{0,1\}~\forall l,t$. When $a_t^l=0$, no repair work is to be carried out at $l$. Similarly, when $a_t^l=1$, repair work is carried out at $l$.
  \item \emph{Simulator T:} The repair time associated with each damaged location depends on the state of the damage to the component at that location (see Table~\ref{T1}). This repair time is random and is exponentially distributed with expected repair times shown in Table~\ref{T1}. Given $s_t$ and $a_t$, $T$ gives us the new state $s_{t+1}$. We say that a repair action is complete as soon as at least one of the locations where repair work is carried out is fully repaired. Let's denote this completion time at every $t$ by $\hat t$. Note that it is possible for the repair work at two or more damaged locations to be completed simultaneously. Once the repair action is complete, the units of resources at remaining locations, where repair work was not complete, are also available for reassignment along with unit of resources where repair was complete. The new repair time at such unrepaired locations is calculated by subtracting $\hat t$ from the time required to repair these locations. It is also possible to reassign the unit of resource at the same unrepaired location if it is deemed important for the repair work to be continued at that location by the planner. Because of this reason, preemption of repair work during reassignment is not a restrictive assumption, on the contrary, it allows greater flexibility to the decision maker for planning.  Because the repair times are random, the outcomes of repair actions are random as not the same damaged component will be repaired first every time (random repair time), even if the same repair action $a_t$ is applied in $s_t$. Hence, our simulator $T$ is stochastic. Alternative formulation where outcome of repair action is deterministic is also an active area of research \cite{our,emi,iEMSs}.
  \item \emph{Rewards R:} We wish to optimally plan decisions so that maximum people will get water in minimum amount of time. We combine these two competing objectives to define our reward as:
       \begin{equation}\label{rew}
         R(s_t,a_t,s_{t+1})=\frac{r}{t_{rep}},
       \end{equation}
        where $r$ is the number of people who have water after action $a_t$ is completed, and $t_{rep}$ is the total repair time (days) required to reach $s_{t+1}$ from any initial state $s_0$. Note that our reward function is stochastic because the outcome of our action $a_t$ is random.
  \item \emph{Initial State I:} We have already described the stochastic damage model of the components for the modeled network in Section~\ref{test2} and Section~\ref{test3}. The initial damage states associated with the components will be provided by these models.
  \item \emph{Discount factor $\gamma$:} In our simulation studies, $\gamma$ is fixed at 0.99.
\end{itemize}

\subsection{Rollout}
The rollout algorithm was first proposed for stochastic scheduling problems by Bertsekas and Castanon \cite{Bertsekas1999}. Instead of the dynamic programming formalism \cite{Bertsekas1999}, we motivate the rollout algorithm in relation to the simulation-based representation of our MDP.
Suppose that we have access to a non-optimal policy $\pi$, and our aim is to compute an improved policy $\pi'$. Then, we have:
\begin{equation}\label{roll}
  \pi'(s_t)=\arg\max_{a_t}Q_\pi(s_t,a_t),
\end{equation}
where the $Q$ function is as defined in \eqref{Qval}. If the policy defined in \eqref{roll} $\pi'$ is non-optimal, it is a \emph{strict} improvement over $\pi$ \cite{howard}. This result is termed as \emph{policy improvement theorem}. Note that the improved policy $\pi'$ is generated as a greedy policy w.r.t. $Q_\pi$. Unlike the exact solution methods described in Section~\ref{sol}, we are interested here in computing $\pi'$ only for the current state. Methods that use \eqref{roll} as the basis for updating the policy suffer from the \emph{curse of dimensionality}. Before performing the policy improvement step in \eqref{roll}, we have to first calculate the value of $Q_\pi$. Calculating the value of $Q_\pi$ in \eqref{roll} is known as \emph{policy evaluation}. Policy evaluation is intractable for large or continuous state and action spaces. Approximation techniques alleviate this problem by calculating an approximate $Q$ value function. Rollout is one such approximation technique that utilises monte-carlo simulations. Particularly, rollout can be formulated as an approximate policy iteration algorithm \cite{Fern,lagoudakis2003reinforcement}. An implementable (programming sense) stochastic function (simulator) $SimQ(s_t,a_t,\pi,h)$ is defined in such a way that its expected value is $Q_\pi(s_t,a_t,h)$, where $h$ is a finite number representing horizon length. In the rollout algorithm, $SimQ$ is implemented by simulating action $a_t$ in state $s_t$ and following $\pi$ thereafter for $h-1$ steps. This is done for all the actions $a_t \in A(s_t)$. A finite horizon approximation of $Q_\pi(s_t,a_t)$ (termed as $Q_\pi(s_t,a_t,h)$), is required; our simulation would never finish in the infinite horizon case because we would have to follow policy $\pi$ indefinitely. However, $V^\pi(s_t)$, and consequently $Q_\pi(s_t,a_t)$, is defined over the infinite horizon. It is easy to show the following:
\begin{equation}\label{Qapprox}
  \left | Q_\pi(s_t,a_t)- Q_\pi(s_t,a_t,h)\right |=\frac{\gamma^{\,h}R_{max}}{1-\gamma}.
\end{equation}
The approximation error in \eqref{Qapprox} reduces exponentially fast as $h$ grows. Therefore, the $h$-horizon results apply to the infinite horizon setting, for we can always choose $h$ such that the error in \eqref{Qapprox} is negligible. To summarize, the rollout algorithm can be presented in the following fashion for our problem:
\begin{algorithm}
    \caption{Uniform Rollout ($\pi$,h,$\alpha$,$s_t$)}
    \label{rollout}
    \begin{algorithmic}
    \For{$i=1$ to $n$}
    \For{$j=1$ to $\alpha$}
    \State $\tilde a^{i,j} \gets SimQ(s_t,a^{i,j}_t,\pi,h)$ \Comment{See algorithm 2}
    \EndFor
    \EndFor
    \State  $\tilde a^i_t \gets \emph{Mean}(\tilde a^{i,j}$)
    \State  $k \gets \arg\max\tilde a^i_t$
    \State \Return $a^k_t$
    \end{algorithmic}
    \end{algorithm}
\begin{algorithm}
    \caption{Simulator $SimQ(s_t,a^{i,j}_t,\pi,h)$}
    \label{sim}
    \begin{algorithmic}
    \State $s_{t+1}\gets T(s_t,a^{i,j}_t)$
    \State $r \gets R(s_t,a^{i,j}_t,s_{t+1})$
    \For{$p=1$ to $h-1$}
    \State $s_{t+1+p} \gets T(s_{t+p},\pi(s_{t+p}))$
    \State $r\gets r+ \gamma^{\,p}R(s_{t+p},\pi(s_{t+p}),s_{t+1+p})$
    \EndFor
    \State \Return $r$
    \end{algorithmic}
    \end{algorithm}

In Algorithm~\ref{rollout}, $n$ denotes $|A(s_t)|$. Note that Algorithm~\ref{sim} returns the discounted sum of rewards. When $h = t_c$, we term the rollout as complete rollout, and when $h < t_c$, the rollout is called truncated rollout \cite{Bertsekas1999}. It is possible to analyse the performance of uniform rollout in terms of uniform allocation $\alpha$ and horizon depth $h$ \cite{Fern,dimitri2008a}.

\subsection{Optimal Computing Budget Allocation}
In the previous section, we presented the rollout method for solving our MDP problem. In the case of uniform rollout, we allocate a fixed rollout sampling budget $\alpha$ to each action, i.e., we obtain $\alpha$ number of rollout samples per candidate action to estimate the $Q$ value associated with the action. In the simulation optimization community, this is analogous to total equal allocation (TEA)\cite{fu} with a fixed budget $\alpha$ for each simulation experiment (a single simulation experiment is equivalent to one rollout sample). In practice, we are only interested in the best possible action, and we would like to direct our search towards the most promising candidates. Also, for large real-world problems, the simulation budget available is insufficient to allocate $\alpha$ number of rollout samples per action. We would like to get a rough estimate of the performance of each action and spend the remaining simulation budget in refining the accuracy of the best estimates. This is the classic exploration vs. exploitation  problem faced in optimal learning and simulation optimization problems.

Instead of a uniform allocation $\alpha$ for each action, non-uniform allocation methods have been explored in the literature pertaining to Algorithm~\ref{rollout} called as \emph{adaptive rollout} \cite{Dimitrakakis2008b}. An analysis of performance guarantees for adaptive rollout remains an active area of research \cite{Dimitrakakis2008b,dimitri2018,lazaric}. These non-uniform allocation methods guarantee performance without a constraint on the budget of rollouts. Hence, we explore an alternative non-uniform allocation method that would not only fuse well into our solutions (adaptively guiding the stochastic search) but would also incorporate the constraint of simulation budget in its allocation procedure. Numerous techniques have been proposed in the simulation optimization community to solve this problem. We draw upon one of the best performers \cite{branke}  that naturally fits into our solution framework---OCBA. Moreover, the probability of correct selection $P\{CS\}$ of an alternative in OCBA mimics finding the best candidate action at each stage in Algorithm~\ref{rollout}. Formally, the OCBA problem \cite{Chen2000} for Section~\ref{probform} can be stated as :
\begin{equation}\label{ocba}
  \max_{N_1,\ldots,N_n}P\{CS\}~ \text{such~that} \sum_{i=1}^{n}N_i=B,
\end{equation}
where $B$ represents the simulation budget for determining optimal $a_t$ for $s_t$ at any $t$, and $N_i$ is the simulation budget for the $i^{th}$ action at a particular $t$. At each OCBA allocation step (for the definition of the allocation step, see variable $l$ in \cite{Chen2000}), barring the best alternative, the OCBA solution assigns an allocation that is directly proportional to the variance of each alternative and inversely proportional to the squared difference between the mean of that alternative and the best alternative.

Here, we only provide information required to initialize the OCBA algorithm. For a detailed description of OCBA, including the solution to the problem in \eqref{ocba}, see \cite{Chen2000}. The key initialization variables, for the OCBA algorithm \cite{Chen2000}, are $k$, $T$ (not to be confused with $T$ in this paper), $\Delta$, and $n_0$. The variable $k$ is equal to variable $n$ in our problem. The value of $n$ changes at each $t$ and depends on the number of damaged components and units of resources. The variable $T$ is equal to per-stage budget $B$ in our problem. More information about the exact value assigned to $B$ is described in Section~\ref{simres}. We follow the guidelines specified in \cite{Chen1999} to select $n_0$ and $\Delta$; $n_0$ in the OCBA algorithm is selected equal to 5, and $\Delta$ is kept at $15\%$ of $n$ (within rounding).
\section{SIMULATION RESULTS}\label{simres}
\begin{figure}
	\centering
	\includegraphics[width=\columnwidth]{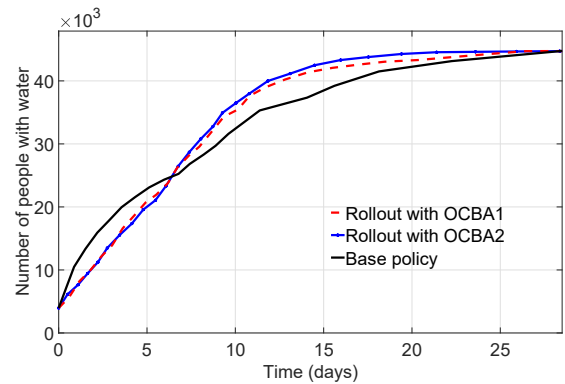}
	\caption{Performance comparison of rollout vs base policy for 3 units of resources.}
	\label{fig4}
\end{figure}

\begin{figure}	
	\centering
	\includegraphics[width=\columnwidth]{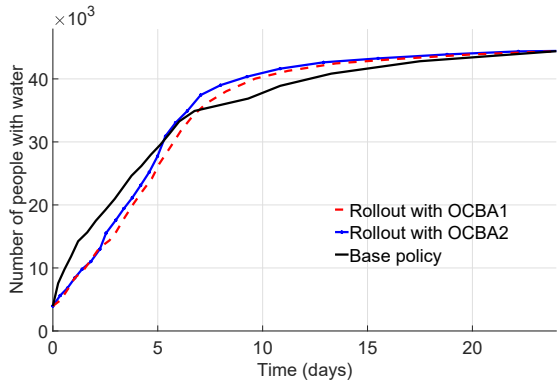}
	\caption{Performance comparison of rollout vs base policy for 5 unit of resources.}
	\label{fig5}
\end{figure}

We simulate 100 different initial damage scenarios for each of the plots presented in this section. There will be a distinct recovery path for each of the initial damage scenarios. All the plots presented here represent the average of 100 such recovery paths. Two different simulation plots of rollout fused with OCBA are provided in Fig.~\ref{fig4} and Fig.~\ref{fig5}. They are termed as rollout with OCBA1 and rollout with OCBA2. The method applied is the same for both cases; only the per-stage simulation budget is different. A per-stage budget (budget at each decision time $t$) of $B=5 \cdot n +5000$ is assigned for rollout with OCBA1 and $B=5 \cdot n +10000$ for rollout with OCBA2.
Fig.~\ref{fig4} compares the performance of rollout fused with OCBA and base policy. The rollout algorithm is known to have the ``lookahead property" \cite{Bertsekas1999}. This behavior of the rollout algorithm is evident in the results in Fig.~\ref{fig4}, where the base policy initially outperforms the rollout policy, but after about six days the former steadily outperforms the later. Recall, that our objective is to perform repair actions so that maximum people will have water in minimum amount of time. Evaluating the performance of our method in meeting this objective is equivalent to checking the area under the curve of our plots. This area represents the product of the number of people who have water and the number of days for which they have water. A larger area represents that greater number of people were benefitted as a result of the recovery actions. The area under the curve for recovery with rollout (blue and red plots) is more than its base counterpart (black). A per-stage budget increase of 5000 simulations in rollout with OCBA2 with respect to rollout with OCBA1 shows improvements in the recovery process.

In the plots shown in Fig.~\ref{fig5}, we use $M=5$. In the initial phase of planning, it might appear that the base policy outperforms the rollout for a substantial amount of time.  However, this is not the case. Note that the number of days for which the base policy outperforms rollout, in both Fig.~\ref{fig4} and Fig.~\ref{fig5}, is about six days, but because the number of resource units has increased from three to five, the recovery is faster, giving an illusion that the base policy outperforms rollout for a longer duration. It was verified that the area under the curve for recovery with rollout (blue and red curves) is more than its base counterpart (black curve). Because OCBA is fused with rollout here, we would like to ascertain the exact contribution of the OCBA approach in enhancing the rollout performance.
\begin{figure}
	\centering
	\includegraphics[width=\columnwidth]{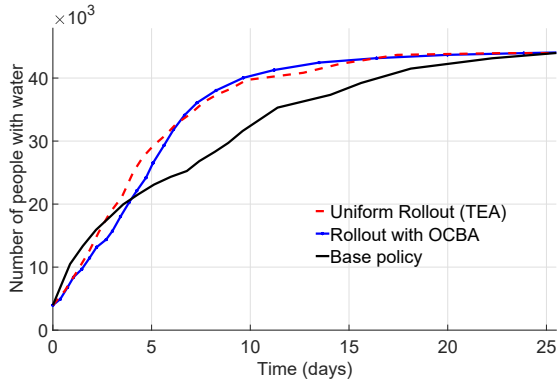}
	\caption{Performance comparison of uniform rollout (TEA), rollout with OCBA and base policy for 3 units of resources.}
	\label{fig6}
\end{figure}

For the rollout with OCBA in Fig.~\ref{fig6}, $B=5 \cdot n +20000$, whereas $\alpha=200$ for the uniform rollout simulations. The recovery as a result of these algorithms outperforms the base policy recovery in all cases. Also, rollout with OCBA performs competitively with respect to uniform rollout despite a meagre simulation budget of $10\%$ of uniform rollout. The area under the recovery process in Fig.~\ref{fig6}, as a result of uniform rollout, is only marginally greater than that due to rollout with OCBA. Note that after six days, OCBA slightly outperforms uniform rollout because it prioritizes the simulation budget on the most promising actions per-stage. Rollout exploits this behavior in each stage and gives a set of sequential recovery decisions that further enhances the outcome of the recovery decisions. We would like to once again stress that such an improvement is being achieved at a significantly low simulation budget with respect to uniform rollout. Therefore, these two algorithms form a powerful combination together, where each algorithm consistently and sequentially reinforces the performance of the other. Such synergistic behavior of the combined approach is appealing. Lastly, our simulation studies show that increments in the simulation budget of rollout results in marginal performance improvement for each increment. Beyond a certain increment in the simulation budget, the gain in performance might not scale with the simulation budget expended. A possible explanation is that small simulation budget increase might not dramatically change the approximation of $Q$ value function associated with a state-action pair. Thus, $\pi'$ in \eqref{roll} might not show a drastic improvement compared to the one computed by a lower simulation budget (policy improvement based on $Q$ approximation that utilises lower simulation budget).

\section{Future Work}
For future work, we would like to leverage the availability of multiple base polices in the aftermath of hazards in our framework and incorporate \emph{parallel rollout} in the solution method \cite{Chang2004}. We anticipate further improvements to the performance demonstrated here when OCBA is fused with parallel rollout. In the future, we will also present the inter-relationship in other critical infrastructure systems like electrical power, roads, bridges, and water networks and the impact such dynamic interactive system has on the recovery process post-hazard. We are also interested in exploring the social impact of the optimized recovery process. We will examine how to incorporate meta-heuristics to guide the stochastic search that determines most promising actions \cite{kaveh}.

\bibliography{IEEEabrv,IEEEexample}

\begin{thebibliography}{10}
\providecommand{\url}[1]{#1}
\csname url@rmstyle\endcsname
\providecommand{\newblock}{\relax}
\providecommand{\bibinfo}[2]{#2}
\providecommand\BIBentrySTDinterwordspacing{\spaceskip=0pt\relax}
\providecommand\BIBentryALTinterwordstretchfactor{4}
\providecommand\BIBentryALTinterwordspacing{\spaceskip=\fontdimen2\font plus
\BIBentryALTinterwordstretchfactor\fontdimen3\font minus
  \fontdimen4\font\relax}
\providecommand\BIBforeignlanguage[2]{{%
\expandafter\ifx\csname l@#1\endcsname\relax
\typeout{** WARNING: IEEEtran.bst: No hyphenation pattern has been}%
\typeout{** loaded for the language `#1'. Using the pattern for}%
\typeout{** the default language instead.}%
\else
\language=\csname l@#1\endcsname
\fi
#2}}

\bibitem{saeed}
S.~Nozhati, B.~R. Ellingwood, H.~Mahmoud, and J.~W. van~de Lindt, ``Identifying
  and {A}nalyzing {I}nterdependent {C}ritical {I}nfrastructure in
  {P}ost-{E}arthquake {U}rban {R}econstruction,'' in \emph{Proc. of the 11th
  Natl. Conf. in Earthq. Eng.}\hskip 1em plus 0.5em minus 0.4em\relax Los
  Angel., CA: Earthq. Eng. Res. Inst., Jun 2018.

\bibitem{ocbar1}
T.~Sun, Q.~Zhao, and P.~B. Luh, ``A {R}ollout {A}lgorithm for {M}ultichain
  {M}arkov {D}ecision {P}rocesses with {A}verage {C}ost,'' in \emph{Posit.
  Syst.}, R.~Bru and S.~Romero-Viv{\'o}, Eds.\hskip 1em plus 0.5em minus
  0.4em\relax Springer Berl. Heidelb., 2009.

\bibitem{ocbarr}
L.~P{\'e}ret and F.~Garcia, ``Online {R}esolution {T}echniques,'' \emph{Markov
  Decis. Process. in Artif. Intell.}, pp. 153--184, 2010.

\bibitem{Gilroy1}
\BIBentryALTinterwordspacing
``The {A}ssociation of {B}ay {A}rea {G}overnments ({City of Gilroy Annex}),''
  2011. [Online]. Available:
  \url{http://resilience.abag.ca.gov/wp-content/documents/2010LHMP/Gilroy-Annex-2011.pdf}
\BIBentrySTDinterwordspacing

\bibitem{Gilroy2}
\BIBentryALTinterwordspacing
{Semseler, R.G. and T. Akel}, ``{2010 Urban Water Management Plan (City of
  Gilroy)}.'' [Online]. Available:
  \url{http://www.ci.gilroy.ca.us/265/Water-Management-Plan}
\BIBentrySTDinterwordspacing

\bibitem{abrahamson}
N.~A. Abrahamson, W.~J. Silva, and R.~Kamai, \emph{Update of the {AS08}
  ground-motion prediction equations based on the {NGA-West2} data set}.\hskip
  1em plus 0.5em minus 0.4em\relax Pac. Earthq. Eng. Res. Cent., 2013.

\bibitem{hazus}
\BIBentryALTinterwordspacing
{Department of Homeland Security, Emergency Preparedness and Response
  Directorate, FEMA, Mitigation Division}, \emph{Multi-hazard Loss Estimation
  Methodology, Earthquake Model: HAZUS-MH MR1, Advanced Engineering Building
  Module}, Wash., DC, Jan 2003. [Online]. Available:
  \url{https://www.hsdl.org/?view&did=11343}
\BIBentrySTDinterwordspacing

\bibitem{adachi}
T.~Adachi and B.~R. Ellingwood, ``Serviceability {A}ssessment of a {M}unicipal
  {W}ater {S}ystem {U}nder {S}patially {C}orrelated {S}eismic {I}ntensities,''
  \emph{Comput.-Aided Civ. and Infrastruct. Eng.}, vol.~24, no.~4, pp.
  237--248, 2009.

\bibitem{ALA}
\BIBentryALTinterwordspacing
J.~Eidinger \emph{et~al.}, ``Seismic fragility formulations for water
  systems,'' \emph{Am. Lifelines Alliance, G\&E Eng. Syst. Inc.}, 2001.
  [Online]. Available: \url{http://homepage. mac. com/eidinger}
\BIBentrySTDinterwordspacing

\bibitem{ouyang}
M.~Ouyang, L.~Due{\~n}as-Osorio, and X.~Min, ``A three-stage resilience
  analysis framework for urban infrastructure systems,'' \emph{Struc. Saf.},
  vol.~36, pp. 23--31, 2012.

\bibitem{masoomi}
H.~Masoomi and J.~W. van~de Lindt, ``Restoration and functionality assessment
  of a community subjected to tornado hazard,'' \emph{Struct. and Infrastruct.
  Eng.}, vol.~14, no.~3, pp. 275--291, 2018.

\bibitem{Puterman}
M.~L. Puterman, \emph{Markov Decision Processes: Discrete Stochastic Dynamic
  Programming}, 1st~ed.\hskip 1em plus 0.5em minus 0.4em\relax N. Y., NY, USA:
  John Wiley \& Sons, Inc., 1994.

\bibitem{howard}
R.~A. Howard, \emph{Dynamic Programming and Markov Processes}.\hskip 1em plus
  0.5em minus 0.4em\relax Camb., MA: MIT Press, 1960.

\bibitem{Bellman}
R.~Bellman, \emph{Dynamic Programming}, 1st~ed.\hskip 1em plus 0.5em minus
  0.4em\relax Princet., NJ, USA: Princet. Univ. Press, 1957.

\bibitem{Shapley}
\BIBentryALTinterwordspacing
L.~S. Shapley, ``Stochastic games,'' \emph{Proc. of the Nat. Acad. of Sci.},
  vol.~39, no.~10, pp. 1095--1100, 1953. [Online]. Available:
  \url{http://www.pnas.org/content/39/10/1095}
\BIBentrySTDinterwordspacing

\bibitem{kallenberg2003finite}
L.~Kallenberg, ``Finite state and action {MDP}s,'' in \emph{Handb. of Markov
  Decis.Process.}\hskip 1em plus 0.5em minus 0.4em\relax Springer, 2003, pp.
  21--87.

\bibitem{Fern}
A.~Fern, S.~Yoon, and R.~Givan, ``Approximate policy iteration with a policy
  language bias: {S}olving relational {M}arkov decision processes,'' \emph{J.
  of Artif. Intell. Res.}, vol.~25, pp. 75--118, 2006.

\bibitem{our}
\BIBentryALTinterwordspacing
S.~Nozhati, Y.~Sarkale, B.~Ellingwood, E.~K.~P. Chong, and H.~Mahmoud,
  ``Near-optimal planning using approximate dynamic programming to enhance
  post-hazard community resilience management,'' \emph{submitted for
  publication}, vol. abs/1803.01451, 2018. [Online]. Available:
  \url{https://arxiv.org/abs/1803.01451}
\BIBentrySTDinterwordspacing

\bibitem{emi}
S.~Nozhati, B.~Ellingwood, H.~Mahmoud, Y.~Sarkale, E.~K.~P. Chong, and
  N.~Rosenheim, ``An approximate dynamic programming approach to community
  recovery management,'' in \emph{Eng. Mech. Inst.}, 2018.

\bibitem{iEMSs}
S.~{Nozhati}, Y.~{Sarkale}, B.~R. {Ellingwood}, E.~K.~P. {Chong}, and
  H.~{Mahmoud}, ``{A combined approximate dynamic programming \& simulated
  annealing optimization method to address community-level food security in the
  aftermath of disasters},'' in \emph{submitted to iEMSs 2018}.

\bibitem{Bertsekas1999}
\BIBentryALTinterwordspacing
D.~P. Bertsekas and D.~A. Castanon, ``Rollout algorithms for stochastic
  scheduling problems,'' \emph{J. of Heuristics}, vol.~5, no.~1, pp. 89--108,
  Apr 1999. [Online]. Available: \url{https://doi.org/10.1023/A:1009634810396}
\BIBentrySTDinterwordspacing

\bibitem{lagoudakis2003reinforcement}
M.~G. Lagoudakis and R.~Parr, ``Reinforcement learning as classification:
  {L}everaging modern classifiers,'' in \emph{Proc. of the 20th Int. Conf. on
  Mach. Learn. (ICML-03)}, 2003, pp. 424--431.

\bibitem{dimitri2008a}
C.~Dimitrakakis and M.~G. Lagoudakis, ``Algorithms and bounds for rollout
  sampling approximate policy iteration,'' in \emph{Recent Adv. in Reinf.
  Learn.}, S.~Girgin, M.~Loth, R.~Munos, P.~Preux, and D.~Ryabko, Eds.\hskip
  1em plus 0.5em minus 0.4em\relax Berl., Heidelb.: Springer Berl. Heidelb.,
  2008, pp. 27--40.

\bibitem{fu}
M.~C. Fu, C.~H. Chen, and L.~Shi, ``Some topics for simulation optimization,''
  in \emph{2008 Winter Simul. Conf.}, Dec 2008, pp. 27--38.

\bibitem{Dimitrakakis2008b}
\BIBentryALTinterwordspacing
C.~Dimitrakakis and M.~G. Lagoudakis, ``Rollout sampling approximate policy
  iteration,'' \emph{Mach. Learn.}, vol.~72, no.~3, pp. 157--171, Sep 2008.
  [Online]. Available: \url{https://doi.org/10.1007/s10994-008-5069-3}
\BIBentrySTDinterwordspacing

\bibitem{dimitri2018}
\BIBentryALTinterwordspacing
------, ``Algorithms and bounds for rollout sampling approximate policy
  iteration,'' \emph{CoRR}, vol. abs/0805.2015, 2008. [Online]. Available:
  \url{http://arxiv.org/abs/0805.2015}
\BIBentrySTDinterwordspacing

\bibitem{lazaric}
\BIBentryALTinterwordspacing
A.~Lazaric, M.~Ghavamzadeh, and R.~Munos, ``Analysis of classification-based
  policy iteration algorithms,'' \emph{J. of Mach. Learn. Res.}, vol.~17,
  no.~19, pp. 1--30, 2016. [Online]. Available:
  \url{http://jmlr.org/papers/v17/10-364.html}
\BIBentrySTDinterwordspacing

\bibitem{branke}
\BIBentryALTinterwordspacing
J.~Branke, S.~E. Chick, and C.~Schmidt, ``Selecting a selection procedure,''
  \emph{Manage. Sci.}, vol.~53, no.~12, pp. 1916--1932, 2007. [Online].
  Available: \url{https://doi.org/10.1287/mnsc.1070.0721}
\BIBentrySTDinterwordspacing

\bibitem{Chen2000}
\BIBentryALTinterwordspacing
C.-H. Chen, J.~Lin, E.~Y{\"u}cesan, and S.~E. Chick, ``Simulation budget
  allocation for further enhancing the efficiency of ordinal optimization,''
  \emph{Discret. Event Dyn. Syst.}, vol.~10, no.~3, pp. 251--270, Jul 2000.
  [Online]. Available: \url{https://doi.org/10.1023/A:1008349927281}
\BIBentrySTDinterwordspacing

\bibitem{Chen1999}
C.-H. Chen, S.~D. Wu, and L.~Dai, ``Ordinal comparison of heuristic algorithms
  using stochastic optimization,'' \emph{IEEE Trans. Robot. Autom.}, vol.~15,
  no.~1, pp. 44--56, Feb 1999.

\bibitem{Chang2004}
\BIBentryALTinterwordspacing
H.~S. Chang, R.~Givan, and E.~K.~P. Chong, ``Parallel rollout for online
  solution of partially observable markov decision processes,'' \emph{Discret.
  Event Dyn. Syst.}, vol.~14, no.~3, pp. 309--341, Jul 2004. [Online].
  Available: \url{https://doi.org/10.1023/B:DISC.0000028199.78776.c4}
\BIBentrySTDinterwordspacing

\bibitem{kaveh}
A.~Kaveh and N.~Soleimani, ``{CBO} and {DPSO} for optimum design of reinforced
  concrete cantilever retaining walls,'' \emph{Asian J. Civ. Eng.}, vol.~16,
  no.~6, pp. 751--774, 2015.

\end{thebibliography}

\end{document}